\documentclass{amsart}

\usepackage{amsmath}
\usepackage{amssymb}
\usepackage{amsthm}
\usepackage{enumerate}
\usepackage{graphics}
\usepackage{pifont}
\usepackage{epsfig} 
\usepackage{setspace} 

\newtheorem{resrep}{Lemma}[section]
\newtheorem{admisimpliesdirich}[resrep]{Proposition}
\newtheorem{measurecounter}[resrep]{Theorem}
\newtheorem{convpoint}[resrep]{Proposition}
\newtheorem{impartto0}[resrep]{Proposition}
\newtheorem{analrep}[resrep]{Proposition}
\newtheorem{admiscounter}[resrep]{Theorem}

\newtheorem{carl}{Lemma}[section]
\newtheorem{normalcounter}[carl]{Theorem}

\newtheorem{carlmult}[carl]{Theorem}

\newtheorem{resrep2}[carl]{Proposition}

\newtheorem{carlBMOA}[carl]{Theorem}
\newtheorem{admisBMOA}[carl]{Proposition}
\newtheorem{limadmis}[carl]{Proposition}
\newtheorem{shiftthm}[carl]{Theorem}

\title[Counterexamples to the weighted Weiss conjecture]{Counterexamples to the discrete and continuous weighted Weiss conjectures}

\begin{document}

\author{Andrew Wynn}
\address{St John's College, Oxford, OX1 3JP}
\email{andrew.wynn@sjc.ox.ac.uk}

\begin{abstract}
Counterexamples are presented to weighted forms of the Weiss conjecture in discrete and continuous time. In particular, for certain ranges of $\alpha$, operators are constructed that satisfy a given resolvent estimate, but fail to be $\alpha$-admissible. For $\alpha \in (-1,0)$ the operators constructed are normal, while for $\alpha \in (0,1)$ the operator is the unilateral shift on the Hardy space $H^2(\mathbb{D})$.
\end{abstract}

\maketitle

\section{Introduction}
Suppose that $(T(t))_{t \geq 0} \subset \mathcal{L}(X)$ is a $C_0$-semigroup with infinitesimal generator $A$ on a Hilbert space $X$. Let $C \in \mathcal{L}(D(A),\mathbb{C})$ be a linear operator which is bounded with respect to the graph norm $\|\cdot \|_{D(A)} := \|A \cdot \|_X + \| \cdot \|_X$ on $D(A)$. Consider the linear system given by
\begin{eqnarray*}
\begin{array}{lllll}
\dot{x}(t) & = & Ax(t), \qquad t > 0; \\
x(0) & = & x_0 \in X;\\
y(t) & = & Cx(t), \qquad t >0.
\end{array}
\end{eqnarray*}
If $x_0 \notin D(A)$, it is not necessarily the case that mild solution $x(t) = T(t)x_0$ lies in $D(A)$ for each $t > 0$ and hence, the \emph{output map} $y(\cdot)$ is not properly defined. However, if it is assumed that $C$ is \emph{admissible} for $A$ in the sense that there exists a constant $M>0$ such that 
\begin{equation*}
\int_0^\infty |CT(t)x_0|^2 dt \leq M^2 \|x_0\|_X^2, \qquad x_0 \in D(A),
\end{equation*}
then the operator $\Psi : D(A) \rightarrow L^2(\mathbb{R}_+)$ given by $(\Psi x)(\cdot):= CT(\cdot)x$ extends continuously to the whole space $X$. In this case, the output map is considered to be given by $y = \Psi x_0$. 

A generalisation of this idea, studied in \cite{HaakLemerdy, HaakKunst}, is to require that $\Psi$ is bounded from $D(A)$ to a weighted $L^2$-space. For $\alpha \in (-1,1)$, the functional $C$ is said to be \emph{$\alpha$-admissible} for $A$ if there exists a constant $M>0$ such that 
\begin{equation} \label{admissibilitydef}
\int_0^\infty t^\alpha |CT(t)x_0|^2 dt \leq M^2 \|x_0\|_X^2, \qquad x_0 \in D(A)
\end{equation}
and it is not difficult to show \cite{HaakLemerdy} that $\alpha$-admissibility implies the resolvent condition
\begin{equation} \label{alpharesolvent}
\sup_{\lambda \in \mathbb{C}_+} (\text{Re} \lambda)^{\frac{1-\alpha}{2}} \|CR(\lambda,A)\|_{X^*} < \infty.
\end{equation}
An interesting problem is to attempt to characterise the operators $A, C$ and weights $\alpha$ for which the reverse implication $(\ref{alpharesolvent}) \Rightarrow (\ref{admissibilitydef})$ is true. The \emph{continuous weighted Weiss conjecture} is said to hold for a class of operators if, given a generator $A$ of that class, $\alpha$-admissibility of any observation operator $C \in \mathcal{L}(D(A),\mathbb{C})$ is equivalent to (\ref{alpharesolvent}).

Initially, the case $\alpha=0$ was considered and in this situation it has been shown that the Weiss conjecture holds whenever $A$ is the generator of a $C_0$-semigroup of contractions \cite{contractions}. However, counterexamples to the unweighted conjecture also exist  \cite{Harper1,oZwart,cZwart}. For a survey of the subject see \cite{Jacobsurvey}. The weighted form of the conjecture was introduced in \cite{HaakLemerdy} for generators of analytic $C_0$-semigroups and for $\alpha \in (-1,1)$ the weighted Weiss conjecture holds in this situation whenever $A^{1/2}$ is admissible for $A$. 

If $A$ is a normal operator generating an analytic $C_0$-semigroup it is easy to check that $A^{1/2}$ is admissible for $A$. Furthermore, if $\alpha \in [0,1)$ and $A$ is a normal operator generating a contractive $C_0$-semigroup, the weighted Weiss conjecture holds without the assumption of analyticity \cite{Wynn}. In \S \ref{sec:setup} it is shown that, even for normal operators, the weighted Weiss conjecture fails in the case $\alpha \in (-1,0)$. 

A discrete form of the weighted Weiss conjecture can also be formulated \cite{Harper2,Wynn}. If $X$ is a Hilbert space and $A \in \mathcal{L}(X)$ with spectrum $\sigma(A) \subset \overline{\mathbb{D}}:=\{ z \in \mathbb{C}: |z| \leq 1\}$ and $C \in X^*$, the linear functional $C$ is said to be discrete $\alpha$-admissible for $A$ if
\begin{equation} \label{discreteadmis}
\sum_{n=0}^\infty (1+n)^\alpha |CA^n x|^2 \leq M\|x\|_X^2, \qquad x \in X
\end{equation}
If $C$ is discrete $\alpha$-admissible for $A$ it can be shown \cite{Wynn} that 
\begin{equation} \label{disalphares}
\sup_{\omega \in \mathbb{D}} (1-|\omega|^2)^{\frac{1-\alpha}{2}} \|C(I-\bar \omega A)^{-1}\|_{X^*} < \infty.
\end{equation}
The \emph{discrete weighted Weiss conjecture} is said to hold for a class of operators if, given a generator $A$ of that class, discrete $\alpha$-admissibility of any $C \in X^*$ is equivalent to (\ref{disalphares}).
If $\alpha=0$, the discrete weighted Weiss conjecture holds for contraction operators \cite{Harper2} and it is shown in \cite{Wynn} that for $\alpha \in (0,1)$ the discrete weighted Weiss conjecture holds for contractive, normal operators. 

In \S \ref{sec:dis:counters} counterexamples are given to the discrete conjecture. It is shown that, even for normal operators, the discrete weighted Weiss conjecture fails for $\alpha \in (-1,0)$. In the case $\alpha \in (0,1)$, the unilateral shift on $H^2(\mathbb{D})$ fails the discrete weighted Weiss conjecture, in contrast to the unweighted case $\alpha =0$. 

\section{Counterexamples to the continuous weighted Weiss conjecture} \label{sec:setup}
Let $\alpha \in (-1,0)$. Suppose that $\mu$ is a finite, positive measure such that $\text{supp}(\mu)$ is a bounded subset of the closed upper half plane $\overline{\Pi_+}:=\{z \in \mathbb{C}: \text{Im}z \geq 0\}$ and $\mu(\mathbb{R})=0$. Let $X := L^2(\Pi_+,\mu)$ and define operators $A \in \mathcal{L}(X), C \in X^*$ by
\begin{equation*}
(Af)(z) :=izf(z), \quad f \in X, z \in \Pi_+, \qquad Cf := \int_{\Pi_+} f(z) d\mu(z), \quad f \in X.
\end{equation*}
Notice that $A \in \mathcal{L}(X)$ is a normal operator, generating a contractive $C_0$-semigroup $(T(t))_{t\geq0}$ on $X$ given by
\begin{equation} \label{sec:counter:semigroupdef}
(T(t)f)(z) = e^{izt} f(z), \qquad f \in X, z \in \Pi_+, t > 0.
\end{equation}
For an interval $I \subset \mathbb{R}$ define $R(I):= \{ x+iy \in \Pi_+: x \in I, y \in (0,|I|/2) \}$. The resolvent estimate (\ref{alpharesolvent}) can be characterised in terms of a bound on $\mu$ on the sets $R(I)$.

\begin{resrep} \label{resrep}
Let $\alpha \in (-1,0)$ and assume that $X, A, C$ and $\mu$ are as above. Then $(\ref{alpharesolvent})$ holds if and only if there exists a constant $c>0$ such that $\mu(R(I)) \leq c|I|^{1+\alpha},$ for any interval $I \subset \mathbb{R}$. 
\end{resrep}
\begin{proof}
For any $\lambda \in \mathbb{C}_+$ and $x \in X$, $CR(\lambda,A)x= \langle R(\lambda,A)x, 1 \rangle_X$ where $1(z)=1, z \in \Pi_+$. Hence, 
\begin{equation*}
\|CR(\lambda,A)\|_{X^*}^2=\|R(\lambda,A)^* 1\|_{X}^2 = \int_{\Pi_+} \frac{d\mu(z)}{|\lambda-iz|^2}, \qquad \lambda \in \mathbb{C}_+.
\end{equation*}
The result follows from \cite{Wynn}, Lemma 5.8. 
\end{proof}

The following result provides a condition on $\mu$ that is necessary for weighted admissibility.

\begin{admisimpliesdirich} \label{admisimpliesdirich}
Let $\alpha \in (-1,0)$ and suppose that $A, C, X$ and $\mu$ are as above. If $C$ is $\alpha$-admissible for $A$ there exists a constant $M>0$ such that 
\begin{equation*}
\left( \int_{\Pi_+} \left| \int_0^\infty e^{izt} t^{\alpha/2} v(t) dt \right|^2 d\mu(z) \right)^{\frac{1}{2}} \leq M \|v\|_{L^2(\mathbb{R}_+)}, \qquad v \in L^2(\mathbb{R}_+).
\end{equation*}
\end{admisimpliesdirich}
\begin{proof}
Suppose that $C$ is $\alpha$-admissible for $A$. Let $v \in L^2(\mathbb{R}_+)$ and 
\begin{equation*}
x \in \mathcal{D}:= \left\{ y \in X : \int_{\Pi_+} \frac{|y(z)|^2}{(\text{Im} z)^{1+\alpha}} d\mu(z) < \infty \right\}.
\end{equation*} 
Then,
\begin{eqnarray*}
\int_{\Pi_+} \int_0^\infty \left| e^{izt} x(z) \right| t^{\frac{\alpha}{2}} |v(t)| dt d\mu(z) &=& \int_{\Pi_+} \left( \int_0^\infty e^{-t \text{Im}z} t^{\frac{\alpha}{2}} |v(t)| dt \right) |x(z)| d\mu(z) \\
\text{(by Cauchy-Schwarz)}&\leq& c_\alpha \|v\|_2 \int_{\Pi_+} \frac{|x(z)|}{(\text{Im}z)^{\frac{1+\alpha}{2}}} d\mu(z)\\
\text{(by Cauchy-Schwarz)} &\leq& c_\alpha \|v\|_2 \left( \mu(\Pi_+) \int_{\Pi_+} \frac{|x(z)|^2}{ (\text{Im}z)^{1+\alpha}} d\mu(z) \right)^{\frac{1}{2}}\\
&<& \infty
\end{eqnarray*}
and Fubini's theorem may be applied. Now,
\begin{eqnarray*}
\left| \left\langle  \left( \int_0^\infty e^{i(\cdot)t} t^{\alpha/2} v(t) dt \right), \bar x(\cdot) \right\rangle_X \right| & = & \left| \int_{\Pi_+} \left( \int_0^\infty e^{izt} t^{\alpha/2} v(t) dt \right) x(z) d\mu(z) \right| \\
\text{(by Fubini)} & = & \left| \left( \int_0^\infty \int_{\Pi_+} e^{izt} x(z) d\mu(z) \right) t^{\alpha/2} v(t) dt \right|\\
(\text{by (\ref{sec:counter:semigroupdef})})& = & \left| \int_0^\infty \big(CT(t)x \big) t^{\alpha/2} v(t) dt \right|\\
\text{(by Cauchy-Schwarz)}&\leq & \left( \int_0^\infty t^\alpha |CT(t)x|^2 dt \right)^{1/2} \|v\|_{L^2(\mathbb{R}_+)}\\
\text{(by $\alpha$-admissibility)} &\leq & M \|x\|_X \|v\|_{L^2(\mathbb{R}_+)}.
\end{eqnarray*}
Since $\mathcal{D}$ is dense in $X$,
\begin{equation*}
\left\| \int_0^\infty e^{i(\cdot)t} t^{\alpha/2} v(t) dt  \right\|_X \leq M \|v\|_{L^2(\mathbb{R}_+)}.
\end{equation*} 
\end{proof}

If $\alpha >-1$ then, upon identifying functions that differ by a constant, the weighted Dirichlet space $\mathcal{D}_{1+\alpha}^2 (\Pi_+)$ contains those analytic functions $F:\Pi_+ \rightarrow \mathbb{C}$ for which
\begin{equation*}
\|F\|_{\mathcal{D}^2_{1+\alpha}}^2:=\int_0^\infty \!\!\! \int_{\!-\infty}^\infty y^{1+\alpha} |F'(x + iy)|^2 dxdy < \infty.
\end{equation*}
Furthermore, by \cite[Theorem 3]{Duren2}, $F \in \mathcal{D}^2_{1+\alpha}(\Pi_+)$ if and only if there exists a function \[w \in L^2(\mathbb{R}_+, t^{-\alpha} dt):= \left\{f:\mathbb{R}_+ \rightarrow \mathbb{C}:f \; \text{measurable}, \; \int_0^\infty t^{-\alpha} |f(t)|^2 dt < \infty \right\}\] and a constant $c \in \mathbb{C}$ with 
\begin{equation} \label{durrep}
 F(z) = \int_0^\infty e^{izt} w(t) dt +c, \qquad z \in \Pi_+.
\end{equation}
In this case there exists a constant $k>0$ with $\|F\|_{\mathcal{D}^2_{1+\alpha}} = k \|w\|_{L^2(\mathbb{R}_+,t^{-\alpha}dt)}.$ Proposition \ref{admisimpliesdirich} now states that for $\alpha \in (-1,0)$, the embedding 
\begin{equation}
\mathcal{D}_{1+\alpha}^2(\Pi_+) \hookrightarrow L^2(\Pi_+,d\mu) \label{diremb}
\end{equation} 
is necessary for $\alpha$-admissibility of $C$ with respect to $A$. Hence, in order to create a counterexample, it is enough to find a measure $\mu$ satisfying $\mu(R(I)) \leq c|I|^{1+\alpha}$ but for which (\ref{diremb}) does not hold. 

In the unweighted case $\alpha=0$, Lemma \ref{resrep} and Proposition \ref{admisimpliesdirich} are still true. However, since the unweighted Weiss conjecture holds for normal operators, a counterexample cannot be created in this case. Indeed, Proposition \ref{admisimpliesdirich} implies that $H^2(\Pi_+) \hookrightarrow L^2(\Pi_+,d\mu)$ is necessary for $0$-admissibility, but by the Carleson measure theorem (see e.g. \cite{Duren1}), this embedding is equivalent to the bound $\mu(R(I)) \leq c |I|$. By Lemma \ref{resrep}, this bound on $\mu$ is the same as (\ref{alpharesolvent}) with $\alpha =0$. In fact, it is for exactly this reason that the unweighted Weiss conjecture is true for normal operators \cite{WeissPowerful}. 

The reason that counterexamples can be found in the case $\alpha \in (-1,0)$ is that measures satisfying (\ref{diremb}) are characterised by a bound involving the Riesz capacity of certain sets
(see e.g. \cite[Theorem 4.4]{Wu}) but not by a simple condition of the form $\mu(R(I)) \leq c |I|^{1+\alpha}$.   

\subsubsection*{Riesz capacities}

Let $\beta \in (0,1)$. The $\beta$-Riesz capacity of a subset $A \subset \mathbb{R}$ is given by 
\begin{equation} \label{sec:counter:capacity}
Cap_\beta (A) := \inf\{\; \|g\|_{L^2(\mathbb{R})}^2 : g \in L^2(\mathbb{R}), I_\beta \ast g \geq 1 \; \text{on} \; A, \; g \geq 0 \; \}
\end{equation} 
where the {\em Riesz kernel} $I_\beta$ is defined by $I_\beta := |x|^{\beta-1}$ (see e.g. \cite{Adams}, p.8). If $O \subset \mathbb{R}$ is an open set define 
\begin{equation*}
R(O):= \bigcup_{i=1}^\infty R(I_i),
\end{equation*}
where $O = \bigcup_{i=1}^\infty I_i$ is the decomposition of $O$ into disjoint open intervals of $\mathbb{R}$. 
If $\alpha \in (-1,0)$, it is shown in \cite[Theorem 6.1]{Dafni}  that there exists a measure $\mu$ on $\Pi_+$ for which $\mu(R(I)) \leq c|I|^{1+\alpha}$ for any interval $I \subset \mathbb{R}$, but for which there does not exist a constant $c>0$ with
\begin{equation} \label{sec:counter:capest}
\mu(R(O)) \leq c \cdot Cap_{-\alpha/2}(O), \qquad O \subset \mathbb{R} \; \; \text{open}.
\end{equation}
Such a measure will be used to construct the counterexamples. With respect to the operators $A$ and $C$ introduced in \S \ref{sec:setup}, it is shown in Lemma \ref{resrep} that the resolvent bound (\ref{alpharesolvent}) is equivalent to the one-box condition $\mu(R(I)) \leq c|I|^{1+\alpha}$, while it will be shown later that the capacity estimate (\ref{sec:counter:capest}) is necessary for $\alpha$-admissibility. However, it will be useful to determine the possible structure of such a measure.

\begin{measurecounter} \label{measurecounter}
Let $\alpha \in (-1,0)$. Then there exists a finite, positive measure $\mu$ on $\Pi_+$ with compact support such that 
\begin{itemize}
\item[(i)] There exists $c>0$ such that $\mu(R(I)) \leq c |I|^{1+\alpha}$, for any interval $I \subset \mathbb{R}$.
\item[(ii)] There does not exist $c>0$ such that $\mu(R(O)) \leq c \cdot Cap_{-\alpha/2}(O)$ holds for every open set $O \subset \mathbb{R}$. 
\end{itemize}
\end{measurecounter}

\begin{proof}
From the proof of \cite[Theorem 6.1]{Dafni} there exists a non-trivial positive measure $\nu$ on $\mathbb{R}$ and a compact set $K \subset \mathbb{R}$ such that: there exists $c>0$ such that $\nu(I) \leq c |I|^{1+\alpha}$, for any interval $I \subset \mathbb{R}$; $Cap_{-\alpha/2}(K) =0$ and $\text{supp}(\nu) \subset K$. 
Since $Cap_{-\alpha/2}(K) =0$, \cite[Theorem 2.3.11]{Adams} implies that there exists a sequence $(O^{(n)})_{n=1}^\infty$ of open sets $O^{(n)} \subset \mathbb{R}$ such that for each $n \in \mathbb{N}$, $O^{(n)} \supseteq O^{(n+1)}, K \subset O^{(n)}$ and additionally,
\begin{equation} \label{capacityrate}
Cap_{-\alpha/2}(O^{(n)}) \leq \frac{1}{n^2}, \qquad n \in \mathbb{N}. 
\end{equation}
The set $O^{(1)}$ can be expressed as a disjoint union of open intervals which form an open cover for $K$. Since $K$ is compact there exists a finite subcover $I_1^{(1)}, \ldots, I_{N_1}^{(1)}$ of intervals such that $\tilde O^{(1)} := \bigcup_{i=1}^{N_1} I_i^{(1)} \supset K$ and since $\tilde O^{(1)} \subset O^{(1)}$ it follows that $Cap_{-\alpha/2}(\tilde O^{(1)}) \leq Cap_{-\alpha/2}( O^{(1)})$. Since $K \subset O^{(2)} \cap \tilde O^{(1)}$, compactness can again be applied and there exist open intervals $I_1^{(2)}, \ldots, I_{N_2}^{(2)}$ such that $\tilde O^{(2)} := \bigcup_{i=1}^{N_2} I_i^{(2)} \supset K$, $\tilde O^{(2)} \subset O^{(2)} \cap \tilde O^{(1)}$ and $Cap_{-\alpha/2}(\tilde O^{(2)}) \leq Cap_{-\alpha/2}( O^{(2)})$. In this way it is possible to inductively define open sets $\tilde O^{(n)} \subset \mathbb{R}$, such that for each $n \in \mathbb{N}$:
\begin{itemize}
\item[(a)] $\tilde O^{(n)} = \bigcup_{i=1}^{N_n} I_i^{(n)}$, for disjoint open intervals $I_1^{(n)}, \ldots, I_{N_n}^{(n)} \subset \mathbb{R}$;
\item[(b)] $K \subset \tilde O^{(n+1)} \subset \tilde O^{(n)} \subset O^{(n)}$;
\item[(c)] $Cap_{-\alpha/2}( \tilde O^{(n)}) \leq Cap_{-\alpha/2}( O^{(n)}) \leq  1/n^2$.
\end{itemize}
For each $n \in \mathbb{N}$, define $\gamma_n := \frac{1}{3} \min_{i=1}^{N_n} |I_i^{(n)}| >0$ and notice that without loss of generality the sets $\tilde O^{(n)}$ can be picked in such a way that $(\gamma_n)_{n=0}^\infty$ is monotone decreasing. Furthermore, since $Cap_{-\alpha/2}(K)=0$ it must be that case that $\gamma_n \rightarrow 0$ as $n \rightarrow \infty$. The measure $\mu$ on $\Pi_+$ is then defined by
\begin{equation} \label{sec:counter:concretedef}
\mu:= \sum_{m=1}^\infty \frac{1}{m^2} \cdot ( \nu \times \delta_{\gamma_m} ),
\end{equation}
where $\delta_x$ is the point mass at $x \in\mathbb{R}$. 

\noindent
(i) Let $I \subset \mathbb{R}$ be an interval. Then,
\begin{align*}
\mu(R(I)) &= \sum_{m=1}^\infty \frac{1}{m^2} \left( \nu \times \delta_{\gamma_m} \right) (R(I))\\
					&\leq \sum_{m=1}^\infty \frac{\nu(I)}{m^2} \leq \frac{c   \pi^2 |I|^{1+\alpha}}{6}.
\end{align*}

\noindent
(ii) For a contradiction suppose that there exists a constant $c>0$ such that $\mu(R(O)) \leq c \cdot Cap_{-\alpha/2}(O)$ for any open set $O \subset \mathbb{R}$. For each $n \in\mathbb{N}$, 
\begin{equation} \label{chopdown}
R(I_{i}^{(n)}) \supset \{ x+iy \in \Pi_+: x \in I_{i}^{(n)}, 0 < y \leq \gamma_n \}, \qquad 1 \leq i \leq N_n
\end{equation} 
and hence,
\begin{eqnarray*}
\mu( R(\tilde O^{(n)})) =	\mu \left( \bigcup_{i=1}^{N_n} R(I_i^{(n)}) \right)
											&=& \sum_{m=1}^\infty \frac{1}{m^2} \left( \nu \times \delta_{\gamma_n} \right) \left( \bigcup_{i=1}^{N_n} R(I_i^{(n)}) \right)\\
(\text{by (\ref{chopdown})})	&\geq& \sum_{m=n}^\infty \frac{1}{m^2} \nu \left( \bigcup_{i=1}^{N_n} I_i^{(n)} \right)\\
(\text{by (a) and (b)}) &\geq& \nu (K) \sum_{m=n}^\infty \frac{1}{m^2} \\
													&\geq& \frac{\nu(K)}{2n}.																					
\end{eqnarray*}													
Hence, from the above inequality and (c),
\begin{equation*}
\frac{\nu(K)}{2n} \leq \mu(R(\tilde O^{(n)})) \leq c \cdot Cap_{-\alpha/2}(\tilde O^{(n)}) \leq \frac{c}{n^2}, \qquad n \in \mathbb{N},
\end{equation*}
contradicting the assumption.
\end{proof}

In view of Proposition \ref{admisimpliesdirich}, to link the measure (\ref{sec:counter:concretedef}) with $\alpha$-admissibility requires linking the capacity estimate (\ref{sec:counter:capest}) with a weighted Dirichlet space. 
In other words, it is useful to relate each function $g \in L^2(\mathbb{R}_+)$ with some $G \in \mathcal{D}_{1+\alpha}^2( \Pi_+)$. To provide this link (see Proposition~\ref{analrep}) it is of interest to derive some properties of the harmonic extension of $I_\beta \ast g$ to the upper half plane $\Pi_+$. The harmonic extension $u_f$ of a function $f \in L^p( \mathbb{R})$ is given by
\begin{equation} \label{harmon}
u_f(x+iy):=(f \ast P_y)(x), \qquad x+iy \in \Pi_+,
\end{equation}
where $P_y(x):= y/\pi(x^2 + y^2)$ is the {\em Poisson kernel} which satisfies
\begin{equation*}
(\mathcal{F}P_y)(t):=  \int_{\!-\infty}^\infty e^{-ist} P_y(s) ds = e^{-|t|y}, \qquad t \in \mathbb{R}, y >0.
\end{equation*}
For suitable functions $g$, the function $Mg$ is the {\em Hardy-Littlewood maximal function} of $g$ defined by
\begin{equation*}
(Mg)(x):= \sup_{I} \frac{1}{|I|} \int_{I} |g(x)| dx, \qquad x \in \mathbb{R}.
\end{equation*}
It is well known (see e.g. \cite{Adams}, p.3) that if $g \in L^p(\mathbb{R})$ for $1<p \leq \infty$ then there exists a constant $c>0$, depending only on $p$, for which 
\begin{equation} \label{hlmaxbdd}
\|Mg\|_{L^p(\mathbb{R})} \leq c \cdot \|g\|_{L^p(\mathbb{R})}.
\end{equation}

\begin{convpoint} [{\cite[Proposition 3.1.2]{Adams}}] \label{convpoint}
Let $\beta \in (0,1)$. Then there exists a constant $A>0$, depending only upon $\beta$ and $p$, such that for any measurable function $g \geq0$ and any $x \in \mathbb{R}$,
\begin{equation*}
(I_{\beta} \ast g)(x) \leq A \|g\|_p^{\beta p} \cdot((Mg)(x))^{1-\beta p}, \qquad 1\leq p < \frac{1}{\beta}.
\end{equation*}
\end{convpoint}

\begin{impartto0} \label{impartto0}
Let $\alpha \in (-1,0)$. Suppose that $g \in L^2(\mathbb{R})$ and let $f:= I_{-\alpha/2} \ast g$. If $u_f$ is the harmonic extension of $f$ to $\Pi_+$, then for any $x \in \mathbb{R}$
\begin{equation} \label{yto0}
|u_f(x+iy)| \rightarrow 0, \qquad y \rightarrow \infty.
\end{equation}
\end{impartto0}
\begin{proof}
Applying Proposition \ref{convpoint} with $\beta=-\alpha/2$ and $p=2$ gives
\begin{equation*}
|f(x)| = |(I_{-\alpha/2} \ast g)(x)| \leq (I_{-\alpha/2} \ast |g|)(x) \leq  A \|g\|_2^{-\alpha} \cdot ((M|g|)(x))^{1+\alpha}, \qquad x \in \mathbb{R}.
\end{equation*}
Hence,
\begin{eqnarray*}
\int_{\!-\infty}^\infty |f(x)|^{\frac{2}{1+\alpha}} dx & = & \int_{\!-\infty}^\infty |(I_{-\alpha/2} \ast g)(x)|^{\frac{2}{1+\alpha}} dx \\
& \leq & A^{\frac{2}{1+\alpha}}  \cdot \|g\|_2^{\frac{-2\alpha}{1+\alpha}}  \int_{\!-\infty}^\infty |(M|g|)(x)|^2 dx\\
(\text{by (\ref{hlmaxbdd})}) & \leq & A^{\frac{2}{1+\alpha}} \cdot \|g\|_2^{\frac{-2\alpha}{1+\alpha}} \cdot c^2 \cdot \|g\|_2^2\\
& = & \tilde A \cdot \|g\|_2^{\frac{2}{1+\alpha}} < \infty.
\end{eqnarray*}
Therefore, $f \in L^{\frac{2}{1+\alpha}}(\mathbb{R})$, where $\frac{2}{1+\alpha} \in (2,\infty)$. It is shown in (\cite{Garnett}, p.17) that $u_f$ must then satisfy 
\begin{equation*}
|u_f(x+iy)| \leq \left(\frac{2}{\pi y} \right)^{\frac{2+\alpha}{2}} \sup_{\eta >0} \left( \int_0^\infty |u_f(x+i\eta)|^{\frac{2}{2+\alpha}}dx \right)^{\frac{2+\alpha}{2}} \rightarrow 0, \qquad y \rightarrow \infty. 
\end{equation*} 
\end{proof}

\begin{analrep} \label{analrep}
Let $\alpha \in (-1,0)$. Suppose that $g \in L^2(\mathbb{R})$, $f:= I_{-\alpha/2} \ast g$ and let $u_f$ be the harmonic extension of $f$ to $\Pi_+$. Then there exists an analytic function $G\in \mathcal{D}_{1+\alpha}^2(\Pi_+)$ and a function $w \in L^2(\mathbb{R}_+,t^{-\alpha} dt)$ for which $\mathrm{Re}G = u_f, \|G\|_{\mathcal{D}_{1+\alpha}^2} = c\|g\|_{L^2(\mathbb{R})}$ and 
\begin{equation}
G(z) = \int_0^\infty e^{izt} w(t) dt, \qquad z \in \Pi_+.
\end{equation}
\end{analrep}
\begin{proof} 
Since $u_f$ is harmonic in $\Pi_+$ there exists an analytic function $\tilde G:\Pi_+ \rightarrow \mathbb{C}$ with $\text{Re}(\tilde G) = u_f$. It is shown in (\cite{Stein}, p.83)  that for any $x+iy \in \Pi_+$,
\begin{equation*}
\int_{\!-\infty}^\infty |\tilde G'(x+iy)|^2 dx = \int_{\!-\infty}^\infty |\nabla u_f(x+iy)|^2 dx = \frac{1}{\pi} \int_{\!-\infty}^\infty |t|^2 |(\mathcal{F}f)(t)|^2 e^{-2\pi |t| y} dt.
\end{equation*}
Furthermore, it is shown in \cite{Okikiolu} that $(\mathcal{F}f)(t) = \cot{(\frac{-\pi\alpha}{4})} |t|^{\alpha/2} (\mathcal{F}g)(t)$ for almost every $t \in\mathbb{R}$. An application of Fubini's theorem implies that
\begin{eqnarray}
\int_0^\infty \!\!\! \int_{\!-\infty}^\infty y^{1+\alpha} |\tilde G'(x+iy)|^2 dx dy & = &  c \int_{\!-\infty}^\infty |t|^{-\alpha} |(\mathcal{F}f)(t)|^2 dt \nonumber \\
& = &  c \int_{\!-\infty}^\infty |(\mathcal{F}g)(t)|^2 dt  \nonumber\\
(\text{by Plancherel}) & = &  2\pi c \|g\|_{L^2(\mathbb{R})}^2 < \infty. \label{Gnorm}
\end{eqnarray}
Hence, $\tilde G \in \mathcal{D}_{1+\alpha}^2(\Pi_+)$ and by (\ref{durrep}) there exists a function $w \in L^2(\mathbb{R}_+,t^{-\alpha}dt)$ and a constant $K \in \mathbb{C}$ for which 
\begin{equation*}
\tilde G(z) = \int_0^\infty e^{izt} w(t) dt + K, \qquad z \in \Pi_+.
\end{equation*}
Furthermore \cite{Duren2}, if $G(z):=\int_0^\infty e^{izt} w(t) dt$, then $G(x+iy)$ has the property that for each $x \in \mathbb{R}$, $G(x+iy) \rightarrow 0, \, y \rightarrow \infty$. By Proposition \ref{impartto0}, $\text{Re}(\tilde G)(x+iy) = u_f(x+iy) \rightarrow 0, y \rightarrow \infty$ and hence, $\text{Re}(K) = \text{Re}(\tilde G - G)(x+iy) \rightarrow 0,\; y \rightarrow \infty$. Since $K$ is constant,
\begin{equation*}
\text{Re}(G) = \text{Re}(\tilde G) = u_f.
\end{equation*}
Finally, since $G' = \tilde G'$, it follows from (\ref{Gnorm}) that $\|G\|_{\mathcal{D}_{1+\alpha}^2} = \sqrt{2\pi c} \cdot \|g\|_{L^2(\mathbb{R})}.$ 
\end{proof}

\subsubsection*{The counterexample} \label{countersect}
It is now possible to prove the main result of this section. Recall that $\alpha \in (-1,0)$, $X:=L^2(\Pi_+,\mu)$, $(Af)(z):=izf(z)$ and $Cf = \int_{\Pi_+} f(z) d\mu(z)$. The argument to show that (\ref{sec:counter:capest}) is necessary for $\alpha$-admissibility of $C$ with respect to $A$ is similar to \cite[Theorem 4.4]{Wu}.

\begin{admiscounter} \label{admiscounter}
Let $\alpha \in (-1,0)$. Suppose that $\mu$ is as in Theorem \ref{measurecounter} and that $X, C$ and $A$ are as above. Then $C$ is not $\alpha$-admissible for $A$ but 
\begin{equation*} 
\sup_{\lambda \in \mathbb{C}_+} (\mathrm{Re}\lambda)^{\frac{1-\alpha}{2}} \|CR(\lambda,A)\|_{X^*} < \infty.
\end{equation*}
\end{admiscounter}

\begin{proof}
Since $\mu$ satisfies property (i) of Theorem \ref{measurecounter}, it follows from Lemma \ref{resrep} that the resolvent estimate (\ref{alpharesolvent}) holds. 

Assume for a contradiction that $C$ is $\alpha$-admissible for $A$ and let $O \subset \mathbb{R}$ be an open set. From the definition of $Cap_{-\alpha/2}(O)$, there exists a function $g \in L^2(\mathbb{R}), g \geq 0$ for which $I_{-\alpha/2} \ast g \geq 1$ on $O$ and $\|g\|_{L^2(\mathbb{R})}^2 \leq 2 Cap_{-\alpha/2}(O)$. Define $f:=I_{-\alpha/2} \ast g$, and let $u_f$ be the harmonic extension of $f$ to the upper half plane. From Proposition \ref{analrep} there exists an analytic function $G:\Pi_+ \rightarrow \mathbb{C}_+$ and a function $w \in L^2(\mathbb{R}_+,t^{-\alpha} dt)$ for which $\text{Re}G= u_f, \|w\|_{L^2(\mathbb{R}_+,t^{-\alpha} dt)} = \|G\|_{D_{1+\alpha}^2} = c\|g\|_{L^2(\mathbb{R})}$ and 
\begin{equation} \label{Grep}
G(z) = \int_0^\infty e^{izt} w(t) dt, \qquad z \in \Pi_+.
\end{equation}
Let $w^\circ(t):=w(t) t^{-\alpha/2}, t \in \mathbb{R}_+$. Then $w^\circ \in L^2(\mathbb{R}_+)$ with $\|w^\circ\|_{L^2(\mathbb{R}_+)} = c\|g\|_{L^2(\mathbb{R})}$ and Proposition \ref{admisimpliesdirich} implies that there exists a constant $M>0$ such that 
\begin{eqnarray}
\int_{\Pi_+} \left| \int_0^\infty e^{izt} t^{\alpha/2} w^\circ(t) dt \right|^2 d\mu(z) & \leq & M^2 \| w^\circ \|_{L^2(\mathbb{R}_+)}^2 \nonumber \\
\Longrightarrow \quad \int_{\Pi_+} \left| \int_0^\infty e^{izt} w(t) dt \right|^2 d\mu(z) & \leq & c^2 M^2 \|g\|_{L^2(\mathbb{R})}^2 \nonumber \\
(\text{by (\ref{Grep})}) \Longrightarrow \quad \int_{\Pi_+} |G(z)|^2 d\mu(z) & \leq & c^2 M^2 \|g\|_{L^2(\mathbb{R})}^2. \label{leqg}
\end{eqnarray}
Now suppose that $O = \bigcup O_j$ where each $O_j$ is an open interval in $\mathbb{R}$. If $x+iy \in R(O_j)$, for $O_j =(a,b)$, then since $f = I_{-\alpha/2} \ast g \geq 1$ on $O_j$,
\begin{equation*}
u_f(x+iy)=(f\ast P_y)(x) \geq (\chi_{O_j} \ast P_y)(x) \geq \int_{x-b}^{x-a} P_y(u) du \geq c\cdot\arctan{(2)}:= \delta>0.
\end{equation*}
This holds for any $x+iy \in R(O_j)$ and hence, 
\begin{equation*}
1 \leq \delta^{-2} |(f \ast P_y)(x)|^2 = \delta^{-2} |u_f(x+iy)|^2, \qquad x+iy \in R(O).
\end{equation*}
Therefore, 
\begin{eqnarray*}
\mu(R(O))  =  \int_{R(O)} d\mu(z)
& \leq & \delta^{-2} \int_{R(O)} |u_f(z)|^2 d\mu(z)\\
(\text{$\text{Re}G = u_f$})& \leq & \delta^{-2} \int_{\Pi_+} |\text{Re} G (z)|^2 d\mu(z)\\
& \leq & \delta^{-2} \int_{\Pi_+} |G(z)|^2 d\mu(z)\\
(\text{by (\ref{leqg})}) & \leq & \delta^{-2} c^2 M^2 \|g\|_{L^2(\mathbb{R})}^2\\
(\text{by definition of $g$}) & \leq & 2 \delta^{-2} c^2 M^2 \cdot Cap_{-\alpha/2}(O).
\end{eqnarray*}
This contradicts property (ii) of Theorem \ref{measurecounter} and hence $C$ is not $\alpha$-admissible for $A$. 
\end{proof}

\section{Counterexamples to the discrete weighted Weiss conjecture} \label{sec:dis:counters}
Suppose that $X$ is a Hilbert space, $A \in \mathcal{L}(X)$ with $\sigma(A) \subset \overline{\mathbb{D}}$ and $C \in X^*$. For $\alpha \in (-1,1)$, it is shown in \cite{Wynn} that if $C$ is discrete $\alpha$-admissible for $A$ then 
\begin{equation} \label{disressec}
\sup_{\omega \in \mathbb{D}} (1-|\omega|^2)^{\frac{1-\alpha}{2}} \|C(I-\bar \omega A)^{-1}\|_{X^*}< \infty.
\end{equation}
If $\alpha \in (0,1)$ and $A$ is a normal operator, it is shown in \cite{Wynn} that $C$ is discrete $\alpha$-admissible for $A$ if and only if (\ref{disressec}) holds. It will be shown that this result fails to generalise in two senses. If $\alpha \in (-1,0)$ there exists a normal operator for which the discrete weighted Weiss conjecture fails. In the case $\alpha=0$, Harper proved in \cite{Harper2} that any contraction operator satisfies the (unweighted) discrete Weiss conjecture. This result fails for $\alpha \in (0,1)$; the unilateral shift on $H^2(\mathbb{D})$, a contractive, non-normal operator, does not satisfy the discrete weighted Weiss conjecture. 

Discrete $\alpha$-admissibility is related to Carleson measures for weighted Dirichlet spaces. For  $\alpha \in (-1,1)$, the weighted Dirichlet space $\mathcal{D}_{\alpha}^2(\mathbb{D})$ contains those analytic functions $f(z)=\sum_{n=0}^\infty f_n z^n$ on $\mathbb{D}$ for which 
\begin{equation*}
\|f\|_{\alpha}^2:=\sum_{n=0}^\infty (1+n)^\alpha |f_n|^2 < \infty. 
\end{equation*}
A positive measure $\mu$ on $\mathbb{D}$ is said to be an $\alpha$-Carleson measure if 
\begin{equation*}
\mathcal{D}_{\alpha}^2(\mathbb{D}) \hookrightarrow L^2(\mathbb{D},\mu):=\{ f : \mathbb{D} \rightarrow \mathbb{C}: f \; \text{measurable}, \; \int_{\mathbb{D}} |f(z)|^2 d\mu(z) < \infty \}.
\end{equation*} 
If $\alpha \in [0,1)$ a measure $\mu$ is a $(-\alpha)$-Carleson measure if and only if there exists a constant $c>0$ such that $\mu(S(I)) \leq c |I|^{1+\alpha}$, for any arc $I \subset \mathbb{T}$ (see, e.g. \cite{Stegenga}). Here, 
\begin{equation*}
S(I) := \{ z=re^{i\theta} :  e^{i\theta} \in I, \; 1- \frac{|I|}{2\pi} \leq r < 1 \}.
\end{equation*}
The following result will be useful in constructing the counterexamples.

\begin{carl} \label{carl}
Let $\alpha \in (-1,1)$. Then there exists a constant $c>0$ such that $\mu(S(I)) \leq c|I|^{1+\alpha}$ for any arc $I \subset \mathbb{T}$ if and only if there exists a constant $k>0$ such that 
\begin{equation} \label{int}
\int_{\mathbb{D}} \frac{d\mu(z)}{|1-\bar \omega z|^2} \leq \frac{k}{(1-|\omega|^2)^{1-\alpha}}, \qquad \omega \in \mathbb{D}.
\end{equation}
\end{carl}
\begin{proof}
It is shown in \cite[\S 3]{Wynn} that (\ref{int}) implies $\mu(S(I)) \leq c |I|^{1+\alpha}$ for any arc $I \subset \mathbb{T}$. For the converse, notice first that by rotational invariance and the fact that $\mu(\mathbb{D}) < \infty$, it is sufficient to show that (\ref{int}) holds for $\omega \in (a,1)$, for some fixed $a \in (0,1)$. Let $\omega > 1/2$. Define arcs $I_n \subset \mathbb{T}$ by $I_n := \{ e^{i\theta} : \theta \in (-2^n\pi (1-\omega), 2^n\pi  (1-\omega) \}$ and sets $A_0:=S(I_0), A_n:=S(I_n) \setminus S(I_{n-1}), n \geq 1$. Notice that for a given $\omega \in (0,1)$, there exists $N_\omega \in \mathbb{N}$ such that $A_n = \emptyset$ for $n \geq N_\omega$. Since $\mu$ satisfies $\mu(S(I)) \leq c|I|^{1+\alpha}$,
\begin{equation} \label{boxest}
\mu(S(I_n)) \leq c |I_n|^{1+\alpha} = \tilde c 2^{n(1+\alpha)} (1-\omega )^{1+\alpha}.
\end{equation}
A simple geometric argument shows that there exists a constant $m>0$, independent of $\omega$, such that 
\begin{equation} \label{distest}
|1-\bar \omega z| \geq \frac{1}{2} |\omega^{-1}-z| \geq \frac{m}{2}  (1-\omega) 2^{n-1}, \qquad z \in A_n, n \geq 1.
\end{equation}
Hence, from (\ref{boxest}) and (\ref{distest}),
\begin{eqnarray*}
\int_{\mathbb{D}} \frac{d\mu(z)}{|1-\bar \omega z|^2} = \sum_{n=0}^\infty \int_{A_n} \frac{d\mu(z)}{|1-\bar \omega z|^2} & \leq & \sum_{n=0}^\infty \frac{4\mu(S(I_n))}{m^2 (1-\omega)^2 2^{2n}} \\
& \leq & \frac{\tilde c}{m^2(1-\omega)^{1-\alpha}} \sum_{n=0}^\infty 2^{n(\alpha-1)}\\
(\alpha \in (-1,1)) & \leq & \frac{k}{(1-\omega^2)^{1-\alpha}}. 
\end{eqnarray*}
\end{proof}

\subsection*{The case $\alpha \in (-1,0)$} 
The idea is the same as in case of continuous weighted admissibility. If $\alpha \in (-1,0)$, there exists a finite, positive Borel measure $\mu$ on $\mathbb{D}$ such that \cite[Theorem 19]{arcozzi}:
\begin{itemize}
\item[(a)] There exists a constant $c>0$ such that $\mu(S(I)) \leq c |I|^{1+\alpha}$, any arc $I \subset \mathbb{T}$;
\item[(b)] $\mathcal{D}_{\!\!-\alpha}^2(\mathbb{D}) \not \hookrightarrow L^2(\mathbb{D},\mu)$.
\end{itemize}
The space $\mathcal{O}(\overline{\mathbb{D}})$ of analytic functions introduced in \cite{Harper2} will be useful. For $\alpha \in (-1,1)$, 
\begin{equation*}
\mathcal{O}(\overline{\mathbb{D}}) := \{ f : \mathbb{D} \rightarrow \mathbb{C} : f \; \text{analytic}, \; \exists R>1 \; \text{with} \; \sum_{n=0}^\infty R^n |f_n| < \infty \}
\end{equation*}
is a dense subspace of $\mathcal{D}_{\alpha}^2(\mathbb{D})$ (\cite{Wynn}, Lemma 2.1). Furthermore, if $A$ is a bounded linear operator on a Hilbert space $X$ with $\sigma(A) \subset \overline{\mathbb{D}}$, then for a function $f(z)= \sum_{n=0}^\infty f_n z^n \in \mathcal{O}(\overline{\mathbb{D}})$ it is possible to define $f(A):= \sum_{n=0}^\infty f_n A^n \in \mathcal{L}(X)$.

\begin{normalcounter}
Let $\alpha \in (-1,0)$. Suppose that $\mu$ is a finite, positive measure on $\mathbb{D}$ satisfying $(a)$ and $(b)$ as above. Let $X:=L^2(\mathbb{D},\mu)$, $(Af)(z):=zf(z), f \in X$ and $Cf:= \int_{\mathbb{D}} f(z) d\mu(z), f \in X$.
Then $C$ is not discrete $\alpha$-admissible for $A$ but
\begin{equation*}
\sup_{\omega \in \mathbb{D}} (1-|\omega|^2)^{\frac{1-\alpha}{2}} \|C(I-\bar \omega A)^{-1}\|_{X^*} < \infty.
\end{equation*}
\end{normalcounter}

\begin{proof}
It is not difficult to show that
\begin{equation} \label{weissharper}
\|Cf(A)\|_{X^*}^2 = \int_\mathbb{D} |f(z)|^2 d\mu(z), \qquad f \in \mathcal{O}(\overline{\mathbb{D}}).
\end{equation}
Since $\mu$ satisfies (a), Lemma \ref{carl} implies that the resolvent estimate holds. 

Suppose for a contradiction that $C$ is discrete $\alpha$-admissible for $A$. Then
there exists a constant $M>0$ such that for any $f = (f_n) \in \ell^2$, 
\begin{eqnarray} \label{ell3}
\left| \sum_{n=0}^\infty (1+n)^{\alpha/2} f_n CA^n x \right|^2  \leq M^2 \|x\|^2 \left(\sum_{n=0}^\infty |f_n|^2\right), \qquad x \in X.
\end{eqnarray}
Suppose that $g(z):= \sum_{n=0}^\infty g_n z^n \in \mathcal{D}_{\!\!-\alpha}^2(\mathbb{D})$ and let $a_n:=(1+n)^{-\alpha/2} g_n, n \in \mathbb{N}$. Then $(a_n)_{n=0}^\infty \in \ell^2$ and (\ref{ell3}) implies that 
\begin{equation*} 
\left| \sum_{n=0}^\infty g_n CA^n x \right|^2 \leq   M^2 \|x\|^2_X \|g\|_{-\alpha}^2, \qquad  x \in X.
\end{equation*}			
Since $\mathcal{O}(\overline{\mathbb{D}}) \subset \mathcal{D}_{\!\!-\alpha}^2(\mathbb{D})$ it follows that $\|Cf(A)\|_{X^*} \leq M \|f\|_{-\alpha}$, for each $f \in \mathcal{O}(\overline{\mathbb{D}})$. From (\ref{weissharper}),
\begin{equation} \label{oleq}
\int_{\mathbb{D}} |f(z)|^2 d\mu(z) \leq M^2 \|f\|_{-\alpha}^2, \qquad f \in \mathcal{O}(\overline{\mathbb{D}}).
\end{equation}
Since $g \in \mathcal{D}_{\!\!-\alpha}^2(\mathbb{D})$, and $\mathcal{O}(\overline{\mathbb{D}})$ is dense in $\mathcal{D}_{\!\!-\alpha}^2(\mathbb{D})$ \cite[Lemma 2.1]{Wynn}, there exist $g^{(n)} \in \mathcal{O}(\overline{\mathbb{D}})$ such that $\|g-g^{(n)}\|_{-\alpha} \rightarrow   0$
as $n \rightarrow \infty$. By Fatou's lemma,
\begin{eqnarray*}
\int_{\mathbb{D}} |g(z)|^2 d\mu(z)  & = & \int_{\mathbb{D}} \liminf_{n \rightarrow \infty} |g^{(n)}(z)|^2 d\mu(z)\\
									& \leq & \liminf_{n \rightarrow \infty} \int_{\mathbb{D}} |g^{(n)}(z)|^2 d\mu(z)\\
									& \leq & M^2 \liminf_{n \rightarrow \infty} \|g^{(n)}\|_{-\alpha}^2 \qquad \text{(by (\ref{oleq}))}\\
									& = & M^2 \|g\|_{-\alpha}^2.
\end{eqnarray*}
Since $g \in \mathcal{D}_{\!\!-\alpha}^2(\mathbb{D})$ was arbitrary, this contradicts the fact that $\mu$ satisfies (b). 
\end{proof}

\subsection*{The case $\alpha \in (0,1)$} 
A simple example of a non-normal contraction operator on a Hilbert space is the unilateral shift $S$ on $H^2(\mathbb{D})$ given by
\begin{equation*}
(Sf)(z):=zf(z), \qquad f \in H^2(\mathbb{D}), z \in \mathbb{D}.
\end{equation*}
Since $S$ is a contraction it satisfies the unweighted discrete Weiss conjecture. However, for $\alpha \in (0,1)$, the resolvent bound
\begin{equation} \label{resollpha}
\sup_{\omega \in \mathbb{D}} (1-|\omega|^2)^{\frac{1-\alpha}{2}} \|C(I-\bar \omega S)^{-1}\| < \infty
\end{equation} 
is not sufficient for discrete $\alpha$-admissibility of an observation operator $C \in H^2(\mathbb{D})^*$ with respect to $S$---see Theorem \ref{shiftthm}. In other words, for $\alpha \in (0,1)$, the discrete weighted Weiss conjecture does not hold for contraction operators. It is possible to translate the counterexample from Theorem \ref{shiftthm} to continuous time operators and deduce that for $\alpha \in (0,1)$, the continuous weighted Weiss conjecture is not true for contractive $C_0$-semigroups. In particular, the right-shift semigroup on $L^2(\mathbb{R}_+)$ does not satisfy the continuous weighted Weiss conjecture for $\alpha \in (0,1)$, which is in contrast to the unweighted case \cite{RightshiftCounter}. This result will be published in a separate paper.

The proof of Theorem \ref{shiftthm} depends upon linking a number of areas of function space theory which are introduced in the following section.

\subsubsection*{Multipliers of Dirichlet spaces, Carleson measures and BMOA}
If $\beta <0$, the Dirichlet space norm $\|\cdot\|_{\beta}^2$ is equivalent \cite{Taylor} to the expression 
\begin{equation} \label{normsequiv}
\int_{\mathbb{D}} |f(z)|^2 (1-|z|^2)^{-(1+\beta)} dA(z), \qquad f \in \mathcal{D}_{\alpha}^2(\mathbb{D})
\end{equation}
where $dA(z) = dxdy, z=x+iy \in \mathbb{D}$ is Lebesgue area measure on $\mathbb{D}$. A function $f$ is said to be a {\em multiplier} from $\mathcal{D}_{\beta}^2(\mathbb{D})$ into $\mathcal{D}_{\gamma}^2(\mathbb{D})$, written $f \in M(\mathcal{D}_{\beta}^2(\mathbb{D}), \mathcal{D}_{\gamma}^2(\mathbb{D}))$, if $fg \in \mathcal{D}_{\gamma}^2(\mathbb{D})$ whenever $g \in \mathcal{D}_{\beta}^2(\mathbb{D})$. Multipliers of Dirichlet spaces are closely related to Carleson measures.
The following result \cite[Theorem 1.1]{Stegenga} is a consequence of the equivalence of (\ref{normsequiv}) to the norm $\|\cdot\|_{\beta}$.

\begin{carlmult} \label{carlmult}
Let $\gamma < \beta \leq 0$. Then $f \in M(\mathcal{D}_{\beta}^2(\mathbb{D}), \mathcal{D}_{\gamma}^2(\mathbb{D}))$ if and only if $f$ is analytic and 
\begin{equation*}
|f(z)|^2 (1-|z|^2)^{-(1+\gamma)} dA(z)
\end{equation*}
is a $\beta$-Carleson measure.
\end{carlmult}

For $\gamma < \beta < 0$ it is shown in \cite{Zhao} that $M(\mathcal{D}_{\beta}^2(\mathbb{D}), \mathcal{D}_{\gamma}^2(\mathbb{D})) = \mathcal{B}^{\frac{2+\beta -\gamma}{2}}(\mathbb{D})$,
where for $\delta>1$, $\mathcal{B}^\delta(\mathbb{D})$ is the {\em weighted Bloch space} of analytic functions $f:\mathbb{D} \rightarrow \mathbb{C}$ for which
\begin{equation*}
\sup_{z \in\mathbb{D}} |f'(z)|(1-|z|^2)^\delta < \infty.
\end{equation*}
The situation is different for multipliers from the Hardy space $\mathcal{D}_{0}^2(\mathbb{D})=H^2(\mathbb{D})$ into a Dirichlet space $\mathcal{D}_{\beta}^2(\mathbb{D})$ for $\beta <0$. In particular, it is shown in \cite{Zhao} that $M(\mathcal{D}_{0}^2(\mathbb{D}), \mathcal{D}_{\beta}^2(\mathbb{D})) = F(2,-\beta,1)$ where the $F$-space $F(p,q,s)$, introduced in \cite{Zhao2}, contains those analytic functions $f:\mathbb{D} \rightarrow \mathbb{C}$ for which 
\begin{equation*}
\sup_{a \in \mathbb{D}} \int_{\mathbb{D}} |f'(z)|^p (1-|z|^2)^q g(z,a)^s dA(z) < \infty.
\end{equation*}
Here $g(z,a)$ is the Green function on $\mathbb{D}$ given by 
\begin{equation*}
g(z,a) := -\log{\left|\frac{a-z}{1-\bar a z}\right|}, \qquad a,z \in \mathbb{D}.
\end{equation*}

In addition to multipliers and Carleson measures, discrete $\alpha$-admissibility with respect to $S$ is related to functions of bounded mean oscillation. For a locally integrable function $f:\mathbb{T} \rightarrow \mathbb{C}$, let $f_I:=\frac{1}{|I|} \int_I f $ denote the mean value of $f$ over the arc $I \subset \mathbb{T}$. Then $f$ is said to have bounded mean oscillation if 
\begin{equation} \label{BMO}
\sup_{I \subset \mathbb{T}} \frac{1}{|I|} \int_I |f(z) - f_I|^2 |dz| < \infty
\end{equation}
and the space $BMOA$ contains those functions in $H^2(\mathbb{D})$ whose boundary functions have bounded mean oscillation. It should be noted (\cite{Zhu2}, p.266) that the space $BMOA$ is unchanged if the $L^2$-norm in (\ref{BMO}) is replaced by an $L^p$-norm, for any $1 \leq p < \infty$. Also, $F(2,0,1)=BMOA$ and for $\beta >0$, the spaces $F(2,\beta,1)$ provide natural generalisations of $BMOA$. The following theorem links $BMOA$ to Carleson measures.

\begin{carlBMOA}[{\cite[Theorem 2]{Jevtic}}] \label{carlBMOA} For $f \in H^2(\mathbb{D})$ the following are equivalent:
\begin{itemize}
\item $f \in BMOA$;
\item For one/all $\beta >0$, the measure $|(\mathcal{I}_\beta f)(z)|^2 (1-|z|^2)^{2\beta-1} dA(z)$ is a $0$-Carleson measure.
\end{itemize}
\end{carlBMOA}
In the above theorem, the fractional derivative operator $\mathcal{I}_\beta: H^2(\mathbb{D}) \rightarrow \mathcal{D}_{\!\!-2\beta}^2(\mathbb{D})$, see (\cite{Zhu1}, p.18), is defined for any $\beta >0$ by 
\begin{equation*}
(\mathcal{I}_\beta f)(z):= \sum_{n=0}^\infty (1+n)^\beta f_n z^n, \qquad f(z):= \sum_{n=0}^\infty f_n z^n \in H^2(\mathbb{D}). 
\end{equation*}
It is also of interest to note that 
\begin{equation} \label{diffs}
(zf(z))' = \sum_{n=0}^\infty (1+n) f_n z^n = (\mathcal{I}_1 f)(z), \qquad f \in H^2(\mathbb{D}), z \in \mathbb{D}.
\end{equation}
It will be shown below that for a linear functional $Cf:= \langle f,c \rangle_{H^2}$, each of the following conditions is equivalent to (\ref{resollpha}) if $\alpha \in (0,1)$:
\begin{itemize}
\item[$\bullet$] $|(\mathcal{I}_1c)(z)|^2 (1-|z|^2) dA(z)$ is a $(-\alpha)$-Carleson measure on $\mathbb{D}$;
\item[$\bullet$] $\mathcal{I}_1 c \in M(\mathcal{D}_{\!\!-\alpha}(\mathbb{D}), \mathcal{D}_{\!\!-2}^2(\mathbb{D})) =\mathcal{B}^{2-\alpha/2}(\mathbb{D})$;
\end{itemize}
and that each of the following conditions is equivalent to discrete $\alpha$-admissibility of $C$ with respect to $S$:
\begin{itemize}
\item[$\bullet$] $\mathcal{I}_{\alpha/2}c \in BMOA$;
\item[$\bullet$] $|(\mathcal{I}_1c)(z)|^2 (1-|z|^2)^{1-\alpha} dA(z)$  is a $0$-Carleson measure on $\mathbb{D}$;
\item[$\bullet$] $\mathcal{I}_1 c \in M(\mathcal{D}_{0}(\mathbb{D}), \mathcal{D}_{\alpha-2}^2(\mathbb{D}))=  F(2,2-\alpha,1)$.
\end{itemize}
For $\alpha \in (0,1)$, it is shown in \cite{Zhao2} that $F(2,2-\alpha,1) \subsetneq \mathcal{B}^{2-\alpha/2}(\mathbb{D})$. 
\subsubsection*{Discrete $\alpha$-admissibility of the unilateral shift}
The first step is to provide an alternative expression for the norm of the operator $C(I-\bar \omega S)^{-1}$. In the following, whenever $C \in H^2(\mathbb{D})$ is an observation operator, $c:=C^*$ is the function in $H^2(\mathbb{D})$ for which $Cf= \langle f,c \rangle_{H^2}, \; f \in H^2(\mathbb{D})$. As pointed out by the referee, Proposition \ref{resrep2} is essentially known; a `folklore' result. I would also like to thank the referee for providing the following short proof.  
\begin{resrep2} \label{resrep2}
Let $C \in H^2(\mathbb{D})^*$. Then for any $\omega \in \mathbb{D}$,
\begin{equation}
\| C(I-\bar \omega S)^{-1}\|_{H^2(\mathbb{D})^*} = \left\| \frac{zc(z) - \omega c(\omega)}{z-\omega} \right\|_{H^2(\mathbb{D})}.
\end{equation}
\end{resrep2}

\begin{proof}
Let $\omega \in \mathbb{D}$, $f \in H^2(\mathbb{D})$ and define $k_\omega(z) := (1-\bar \omega z)^{-1}, z \in \mathbb{D}$. Then 
$C(I-\bar \omega S)^{-1}f = \langle k_\omega f , c \rangle_{H^2}$ and hence, if $P_+:L^2(\mathbb{T}) \rightarrow H^2(\mathbb{D})$ is the Hilbert space orthogonal projection onto $H^2(\mathbb{D})$, it follows that
\begin{equation} \label{resrep12}
C(I-\bar \omega S)^{-1}f = \langle f, \bar k_\omega c \rangle_{L^2} = \langle f, P_+(\bar k_\omega c) \rangle_{H^2}, \qquad f \in H^2(\mathbb{D}).
\end{equation}
Now, $P_+(\bar k_\omega c) = \bar k_\omega c -g$, where $g \in H^2(\mathbb{D})^{\perp} \subset  L^2(\mathbb{T})$ is the unique vector for which $\bar k_\omega c -g \in H^2(\mathbb{D})$. It is easy to check that $g(z):=\frac{\omega c(\omega)}{z-\omega} \in H^2(\mathbb{D})^\perp$ has these properties since
\begin{equation*}
(\bar k_\omega c)(z) = \frac{c(z)}{1-\omega \bar z} = \frac{zc(z)}{z - \omega}, \qquad z \in \mathbb{T}
\end{equation*}
and $z \mapsto \frac{zc(z)-\omega c(\omega)}{z-w} \in H^2(\mathbb{D})$. Therefore,  
\begin{equation*}
P_+(\bar k_\omega c) = \frac{zc(z) - \omega c(\omega)}{z-\omega} 
\end{equation*}
and the result follows from (\ref{resrep12}). 
\end{proof}

\begin{limadmis} \label{limadmis}
Let $\alpha \in (0,1)$ and suppose that $C \in H^2(\mathbb{D})^*$. Then $(\ref{resollpha})$ holds if and only if $\mathcal{I}_1 c \in \mathcal{B}^{2-\frac{\alpha}{2}}(\mathbb{D})$. 
\end{limadmis}
\begin{proof}
Proposition \ref{resrep2} implies that (\ref{resollpha}) holds if and only if 
\begin{equation} \label{resmod}
\sup_{\omega \in\mathbb{D}} (1-|\omega|^2)^{1-\alpha} \int_0^{2\pi} \left| \frac{e^{i\theta}c(e^{i\theta}) - \omega c(\omega)}{e^{i\theta}-\omega} \right|^2 \frac{d\theta}{2\pi} < \infty
\end{equation}
where $c(e^{i\theta}) \in L^2(\mathbb{T})$ is the boundary function of $c \in H^2(\mathbb{D})$. 
It is shown in (\cite{Zhu1}, p.165) that for $f \in H^2(\mathbb{D})$, 
\begin{equation*}
\int_0^{2\pi} \left| \frac{ f(e^{i\theta}) - f(\omega)}{e^{i\theta}- \omega} \right|^2 d\theta \sim \int_\mathbb{D} \frac{|f'(z)|^2 (1-|z|^2)}{|1-\bar \omega z|^2} dA(z), \qquad \omega \in \mathbb{D}.
\end{equation*}
Hence, from the above equivalence, (\ref{diffs}) and (\ref{resmod}), the resolvent bound (\ref{resollpha}) holds if and only if 
\begin{equation*}
\sup_{\omega \in\mathbb{D}} (1-|\omega|^2)^{1-\alpha} \int_\mathbb{D} \frac{|(\mathcal{I}_1c)(z)|^2 (1-|z|^2) dA(z)}{|1-\bar \omega z|^2} < \infty. 
\end{equation*}
Lemma \ref{carl} implies that this holds if and only if $|(\mathcal{I}_1c)(z)|^2 (1-|z|^2) dA(z)$ is a $(-\alpha)$-Carleson measure and by Theorem \ref{carlmult} 
this is equivalent to  
\begin{equation*}
\mathcal{I}_1 c \in M(\mathcal{D}_{\!\!-\alpha}^2(\mathbb{D}), \mathcal{D}_{\!\!-2}^2(\mathbb{D}))=\mathcal{B}^{2-\alpha/2}(\mathbb{D}). 
\end{equation*} 
\end{proof}

\begin{admisBMOA} \label{admisBMOA}
Let $\alpha \in (0,1)$ and suppose that $C \in X^*$ is an observation operator. Then $C$ is discrete $\alpha$-admissible for $A$ if and only if $\mathcal{I}_1 c \in F(2,2-\alpha,1)$.
\end{admisBMOA}

\begin{proof}
Since $F(2,2-\alpha,1)= M(\mathcal{D}_{0}^2(\mathbb{D}), \mathcal{D}_{\alpha-2}^2(\mathbb{D}))$, it follows from Theorem \ref{carlmult} that $\mathcal{I}_1 c \in F(2,2-\alpha,1)$ if and only if 
$|(\mathcal{I}_1 c)(z)|^2 (1-|z|^2)^{1-\alpha} dA(z)$ is a $0$-Carleson measure. Since $\mathcal{I}_{1-\alpha/2}\mathcal{I}_{\alpha/2}= \mathcal{I}_1$, this is the same as saying
\begin{equation*}
|(\mathcal{I}_{1-\alpha/2} ( \mathcal{I}_{\alpha/2}c))(z)|^2 (1-|z|)^{2(1-\alpha/2)-1} dA(z)
\end{equation*}
is a $0$-Carleson measure. By Theorem \ref{carlBMOA} this is equivalent to $\mathcal{I}_{\alpha/2} c \in BMOA$. It is shown in (\cite{Peller}, p.284) that $\mathcal{I}_{\alpha/2} c \in BMOA$ if and only if the generalised Hankel operator $\Gamma_c^\alpha : \ell^2 \rightarrow \ell^2$ represented by the matrix 
\begin{equation*}
\left( \begin{array}{cccc}
c_0 & c_1 & c_2 & \cdots\\
2^{\frac{\alpha}{2}} c_1 & 2^{\frac{\alpha}{2}}c_2 & 2^{\frac{\alpha}{2}} c_3 & \cdots\\
3^{\frac{\alpha}{2}} c_2 & 3^{\frac{\alpha}{2}} c_3 & 3^{\frac{\alpha}{2}} c_4 & \cdots\\
4^{\frac{\alpha}{2}} c_3 & 4^{\frac{\alpha}{2}} c_4 & 4^{\frac{\alpha}{2}} c_5 & \cdots\\
\vdots & \vdots & \vdots &
\end{array} \right)
\end{equation*}
is bounded. Now, if  $f \in H^2(\mathbb{D})$ is given by $f(z):=\sum_{n=0}^\infty f_n z^n, \; z \in \mathbb{D}$, then 
\begin{align*}
\sum_{n=0}^\infty (1+n)^\alpha |CS^nf|^2 &= \sum_{n=0}^\infty (1+n)^{\alpha} \left| \langle S^n f, c \rangle_{H^2} \right|^2\\
&= \sum_{n=0}^\infty \left| \sum_{m=0}^\infty (1+n)^{\alpha/2} f_m \bar c_{n+m}\right|^2\\
&= \|\Gamma_c^\alpha((\bar f_n)_{n=0}^\infty)\|_2^2.
\end{align*}
Hence, boundedness of $\Gamma_c^\alpha$ is equivalent to discrete $\alpha$-admissibility of $C$ with respect to $S$.\end{proof}

\begin{shiftthm} \label{shiftthm} Let $\alpha \in (0,1)$. Then,
\begin{itemize}
\item[(i)] If $(\ref{resollpha})$ holds for an observation operator $C \in X^*$, $C$ is discrete $\beta$-admissible for $S$ for any $\beta \in [0,\alpha)$. 
\item[(ii)] There exists an observation operator $C \in X^*$ which satisfies $(\ref{resollpha})$ but for which $C$ is not discrete $\alpha$-admissible for $S$. 
\end{itemize}
\end{shiftthm}

\begin{proof}
(i) Let $\beta \in [0,\alpha)$. It is shown in \cite{Zhao2} that $\mathcal{B}^{2-\frac{\alpha}{2}} \subset F(2,2-\beta,1)$ and hence, by Propositions \ref{limadmis} and \ref{admisBMOA}, $C$ is discrete $\beta$-admissible for $S$. 

(ii) Since $F(2,2-\alpha,1)  \subsetneq B^{2-\frac{\alpha}{2}}$, there exists a function 
\begin{equation*}
f(z):=\sum_{n=0}^\infty f_n z^n \in B^{2-\frac{\alpha}{2}} \setminus F(2,2-\alpha,1).
\end{equation*}
In particular, $f$ is analytic on $\mathbb{D}$,
\begin{equation*}
\int_{\mathbb{D}} |f(z)|^2 (1-|z|^2) dA(z) \leq k \int_{\mathbb{D}} (1-|z|^2)^{-1+\alpha} dA(z) < \infty
\end{equation*}
and by (\ref{normsequiv}), $f \in \mathcal{D}_{\!\!-2}^2(\mathbb{D})$. Since $\mathcal{I}_1 : H^2(\mathbb{D}) \rightarrow \mathcal{D}_{\!\!-2}^2(\mathbb{D})$ is an isomorphism,
\begin{equation*}c(z):= \sum_{n=0}^\infty \frac{f_n}{1+n} \cdot z^n \in H^2(\mathbb{D})
\end{equation*}
and $\mathcal{I}_1 c = f$. Hence, $Cg:=\langle g,c \rangle_{H^2}, \; g \in H^2(\mathbb{D})$ defines a bounded linear functional on $H^2(\mathbb{D})$. By Proposition \ref{limadmis}, $C$ satisfies (\ref{resollpha}) but by Proposition \ref{admisBMOA}, $C$ is not discrete $\alpha$-admissible for $S$. 
\end{proof}


\end{document}